\documentclass[reqno,12pt]{amsart}

\usepackage{epsf}
\usepackage{graphics}
\usepackage{graphicx}
\usepackage{amssymb}
\usepackage{amsmath}

\date{}

\theoremstyle{plain}
\newtheorem{theorem}{Theorem}

\newtheorem{rem}{Remark}
\newtheorem{question}{Question}
\newenvironment{Rem}{\begin{rem}\rm}{\end{rem}}

\theoremstyle{definition}

\theoremstyle{remark}

\def\N{{\mathbb N}}
\def\Z{{\mathbb Z}}

\def\R{{\mathbb R}}

\newcommand{\cro}[1]{\text{cr}(#1)}

\title{The crossing number of polynomial curve systems}

\author{Sebastian Baader, Jasmin J\"org, Hugo Parlier}
\thanks{This project was supported by ANR-SNF Grant number 200021E\_238147 (SUGAR)}

\begin{document}

\begin{abstract} We determine the crossing number of polynomial size curve systems on standard surfaces, in terms of the genus, up to high precision.

\end{abstract}

\maketitle

\section{Introduction}

The problem of minimising the number of crossings is very old in the field of knot theory, and relatively new in the context of curves on surfaces. Given a closed standard surface $\Sigma_g$ of genus~$g$, and a family $\Gamma=\{\gamma_1,\ldots,\gamma_m\}$ of~$m$ pairwise non-isotopic simple closed curves on it, the crossing number of $\Gamma$ is defined as
$$\text{cr}(\Gamma)=\sum_{k<l} i(\gamma_k,\gamma_l),$$
where $i(\alpha,\beta)$ denotes the geometric intersection number of pairs of curves $\alpha,\beta$ in $\Sigma_g$. 

The following lower bound
$$\text{cr}(\Gamma) > \frac{1}{256(g-1)}\left(m \log \left(\frac{m}{(2g-1)e^6}\right)\right)^2$$
was recently established by Hubard and the third author~\cite{HP}, inspired by crossing numbers of graphs drawn in the plane (see \cite{PTT}). The bound holds for all closed curves, not only for simple ones, and its asymptotic behavior as the number of curves grows is coherent with the exponential growth of Huber's prime geodesic theorem \cite{H}, rather than Mirzakhani's count of simple closed geodesics \cite{M}. Interestingly, this bound predates various families of examples demonstrating its strength, in particular families of arithmetic and sparse curve systems, respectively~\cite{BP,BJK}. Our main result adds to these examples. Considering systems of pairwise non-isotopic simple curves $\Gamma=\{\gamma_1,\ldots,\gamma_m\}$ on $\Sigma_g$ as above, we define
$$\text{Cr}(g,m)=\min\{\text{cr}(\Gamma) \mid \#\Gamma=m\}.$$

\begin{theorem}
\label{polynomial}
For all $\alpha \geq 0$, there exists a constant $N \in \N$, so that for all $g \geq N$
$$\frac{1}{257} \alpha^2 g^{1+2\alpha} (\log(g))^2 \leq \mathrm{Cr}(g,[g^{1+\alpha}]) \leq \frac{9}{4} \alpha^2 g^{1+2\alpha} (\log(g))^2.$$ 
\end{theorem}

The limit case $\alpha=0$ is consistent with the fact that $\Sigma_g$ contains $3g-3$ pairwise non-isotopic disjoint simple curves. More interestingly, the special case $\alpha=\frac{1}{3}$ is consistent with the recent result of Burrin and the third author on the crossing number of systoles on surfaces defined by congruence lattices in $\text{SL}(2,\Z)$~\cite{BP}. In fact, our result took inspiration from the latter. The restriction of Theorem~\ref{polynomial} to curve systems of polynomial growth is essential. Indeed, when the number of curves becomes super-polynomial, the behavior of the crossing number changes abruptly, as described in~\cite{BJK}.

More generally, crossing numbers have different behaviors depending on how the topology relates to the number of curves. As pointed out above, for fixed topology and increasing number of curves the crossing inequality from \cite{HP} is roughly sharp, but this requires looking at all closed curves. It is much more surprising, at least to us, that the inequality remains roughly sharp in other regimes, namely when the number of curves is a function of topology and even under a restriction to simple curves. Indeed, the $\log$ part of the crossing number bound seems destined to be linked to the inverse function of the exponential growth of curves. 

We now proceed to the proof of Theorem~\ref{polynomial} which has two parts, the lower and the upper bounds, presented in the next two sections and in this order.

\section{Bounding the crossing number from below}

Let $\Gamma=\{\gamma_1,\ldots,\gamma_m\}$ be a system of $m$ pairwise non-isotopic simple curves on $\Sigma_g$, with $m=[g^{1+\alpha}]$, where $[-]$ denotes the floor function.
The crossing number bound by Hubard and the third author yields:

\begin{align*}
\cro{\Gamma} &> \frac{1}{256(g-1)}\left([g^{1+\alpha}] \log \left(\frac{[g^{1+\alpha}]}{(2g-1)e^6}\right)\right)^2 \\
&> \frac{1}{256g}\left((g^{1+\alpha}-1) \log \left(\frac{g^{1+\alpha}-1}{(2g-1)e^6}\right)\right)^2.
\end{align*}

We estimate the logarithmic factor:
\begin{align*}
    \log\left(\frac{g^{1+\alpha}-1}{(2g-1)e^6} \right) &=  \log\left(\frac{(2g-1)\frac{1}{2}g^{\alpha}+\frac{1}{2}g^{\alpha}-1}{(2g-1)e^6} \right) \\
    & > \log\left(\frac{(2g-1)\frac{1}{2}g^{\alpha}}{(2g-1)e^6} \right) \\
    &= \log\left(\frac{g^{\alpha}}{2e^6} \right) \\
    &=  \alpha\log(g)-\log\left(2e^6\right),
\end{align*}
where we use that $\frac{1}{2} g^{\alpha}-1 > 0$ for large $g$. From this, we obtain
\begin{align*}
    \cro{\Gamma} &\geq \frac{1}{256} \frac{(g^{1+\alpha}-1)^2}{g}\left( \alpha\log(g)-\log\left(2e^6\right) \right)^2 \\
    &\geq \frac{1}{257} \frac{g^{2+2\alpha}}{g}\alpha^2(\log(g))^2, \\
\end{align*}
which settles the first inequality of Theorem~\ref{polynomial}.

\section{Polynomial families of curves}

\begin{figure}[htb]
\centering
\includegraphics[width=0.5\linewidth]{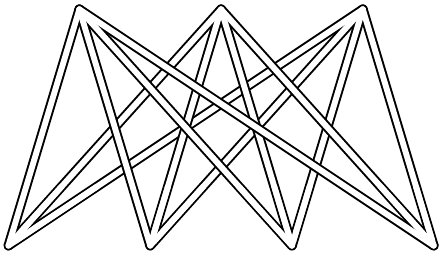}
\caption{Surface $\Sigma(3,4)$}
\label{figure:surface}
\end{figure}

The heart of the construction needed to derive the second inequality of Theorem~\ref{polynomial} is a family of surfaces $\Sigma(p,q) \subset S^3$, the so-called fibre surfaces of torus links of type $T(p,q) \subset S^3$. More precisely, we will need a combinatorial realisation of these surfaces introduced in~\cite{Ba}, also reminiscent of a construction of surfaces with a large number of systoles by the third author~\cite{Pa}. The surface $\Sigma(p,q)$ is a union of $pq$ ribbons located in the neighbourhood of a complete bipartite graph $K_{p,q} \subset \R^3$, whose two vertex sets lie on two skew lines~$U$ and~$L$, $p$ on the `upper' line $U$, and $q$ on the `lower' line $L$. The example depicted in Figure~\ref{figure:surface}, $\Sigma(3,4)$, should be general enough to understand the precise definition of $\Sigma(p,q)$, given in~\cite{Ba}. The concrete embedding of the surface $\Sigma(p,q)$ into $S^3$ is only used to simplify the visualisation of certain curve systems; the final surface $\Sigma_g$ will again be viewed as an abstract surface, rather than an embedded one.

\smallskip

The Euler characteristic of  
$\Sigma(p,q)$ is
$$\chi(\Sigma(p,q))=\chi(K_{p,q})=p+q-pq.$$
For later use, we note that the boundary of $\Sigma(p,q)$ is a torus link of type $T(p,q)$, the number of boundary components of $\Sigma(p,q)$ is $\gcd(p,q)$. In particular, if $\gcd(p,q)=1$, the boundary of $\Sigma(p,q)$ is a knot, i.e.\ it is connected.

Next, we define a system of simple closed curves $\Gamma(p,q)$ on $\Sigma(p,q)$, for all even $q=2k$ with odd $k$. Every member of $\Gamma(p,q)$ is the boundary curve of a subsurface $\Sigma(2,k) \subset \Sigma(p,q)$, with $2$ vertices in a row on the upper line~$U$, and $k$ arbitrary vertices on the lower line~$L$. Since $\gcd(2,k)=1$, the boundary of the subsurface is indeed connected, and the number of curves in $\Gamma(p,q)$ is
$$m(p,q)=(p-1) \binom{2k}{k}.$$
Here, the first factor expresses the number of choices of $2$ consecutive points among the $p$ points on the upper line~$U$, while the second factor expresses the number of choices of $k$ arbitrary points among the $2k$ points on the lower line~$L$. 
To see that all curves in $\Gamma(p,q)$ are homotopically distinct, we distinguish any two curves by their intersection number with a non-trivial arc in $\Sigma(p,q)$. Indeed, for two curves $\gamma, \delta \in \Gamma(p,q)$, there exists a vertex on $U$ or $L$ that is among the defining vertices for $\gamma$ but not for $\delta$. In a neighbourhood of that vertex, we may define a non-trivial simple arc $\alpha$ with endpoints on the boundary of $\Sigma(p,q)$ that intersects $\gamma$ non-trivially. By construction, $\delta$ does not pass through this neighbourhood and is therefore disjoint from $\alpha$.

\begin{figure}[h]
    \centering
    \includegraphics[width=0.5\linewidth]{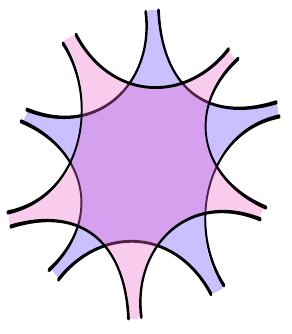}
    \caption{Intersection of two curves in a common vertex}
    \label{figure:intersection}
\end{figure}

In order to estimate the crossing number $\cro{\Gamma(p,q)}$ of the curve system $\Gamma(p,q)$, we make two observations:

\begin{enumerate}
\item Pairs of curves in $\Gamma(p,q)$ with disjoint sets of vertices on the upper line can be deformed to be disjoint: the only points at which the curves potentially intersect lie in common vertices on the lower line; they can easily be removed by a small perturbation.
\item The maximal intersection number between pairs of curves in $\Gamma(p,q)$ is $2k$ or $4k$, depending on whether the curves share one or two vertices on the upper line: even if two curves share certain edges, they can be arranged to be disjoint except in the vertices on the upper line. The local contribution to the intersection number near a vertex on the upper line is not more than $2k$, see Figure~\ref{figure:intersection}.
\end{enumerate}

These two items provide an upper bound for the crossing number of $\Gamma(p,q)$. Let $\gamma \in \Gamma(p,q)$ be a fixed curve. Due to the first item, for a curve to intersect $\gamma$, their vertices on the upper line $U$ need to have a non-trivial overlap, i.e.\ the two curves share one or two vertices on $U$. There are $\binom{2k}{k}$ curves that share two vertices with $\gamma$ on $U$ and at most $2\binom{2k}{k}$ curves that share one vertex with $\gamma$ on $U$. Together with the second item, this yields the following upper bound on the crossing number of the curve system $\Gamma(p,q)$:
$$\cro{\Gamma(p,q)} \leq \left(2 \binom{2k}{k} \cdot 2k+\binom{2k}{k} \cdot 4k \right) \cdot \frac{m(p,q)}{2}=\frac{4k}{p-1}m(p,q)^2.$$

As we will see in the following, a clever choice of $p$ and $q$ settles the second inequality of Theorem~\ref{polynomial}. Recalling that the case $\alpha=0$ corresponds to a simple geometric fact, we fix $\alpha>0$ and $g \in \N$. We would like to set $k= \frac{3}{4} \alpha \log(g)$ and $p=\frac{g}{k}$. For the curve system $\Gamma(p,q)$ to be well-defined, we need $k$ and $p$ to be integers and $k$ to be odd. Therefore, we choose $k$ to be the greatest odd integer smaller than or equal to the desired value. This leads us to define 
\begin{align*}
    k &= 2 \left[\frac{1}{2} \left( \frac{3}{4}\alpha \log(g) -1 \right)\right]+1, \\
    p &= \left[ \frac{4g}{3 \alpha \log(g)} \right]+2, \\
    q &= 2k.
\end{align*}
We will use the following estimates:
\begin{align*}
    \frac{3}{4} \alpha \log(g)-2&\leq k \leq \frac{3}{4} \alpha \log(g),\\
    \frac{4g}{3\alpha \log(g)}+1 &\leq p \leq \frac{4g}{3\alpha \log(g)}+2.
\end{align*}

Our choice implies that the Euler characteristic of the surface $\Sigma(p,q)$ satisfies
\begin{align*}
    |\chi(\Sigma(p,q))| & \leq (p-1)(q-1)\\
    &\leq \left(\frac{4g}{3 \alpha \log(g)} +1\right) \left( \frac{3\alpha \log(g)}{2}-1\right) \\
    &= 2g - \frac{4g}{3 \alpha \log(g)} + \frac{3\alpha \log(g)}{2} -1.
\end{align*}
For large $g$, this is strictly less than $|\chi(\Sigma_g)|=2g-2$. Since we ultimately want a curve system on the closed surface $\Sigma_g$ of genus $g$, we embed the surface $\Sigma(p,q)$ into $\Sigma_g$ in an arbitrary way. 

The above construction thus provides us with a curve system $\Gamma(p,q)$ on $\Sigma_g$ with
$$m(p,q)=(p-1) \binom{2k}{k}$$
curves. The asymptotics of the central binomial coefficient $\binom{2k}{k} \approx \frac{4^k}{\sqrt{\pi k}}$ implies that $m(p,q)$ grows like $\frac{g}{k} \frac{4^k}{\sqrt{\pi k}}$, which is strictly larger than $g^{1+\alpha}$, for all but finitely many $g \in \N$, as we will see in the following.
\begin{align*}
    \binom{2k}{k} &\geq c \frac{4^k}{\sqrt{k}} \geq c \frac{4^{\frac{3}{4}\alpha \log(g)-2}}{\sqrt{\frac{3}{4} \alpha \log(g)}} \\
    &= c'  \frac{4^{\frac{3}{4}\alpha \log(g)}}{\sqrt{ \alpha \log(g)}} = c' \frac{g^{\frac{3}{4}\log(4)\alpha}}{\sqrt{ \alpha \log(g)}},
\end{align*}
where the constants can be taken to be $c=\frac{7}{8 \sqrt{\pi}}$ and $c'= \frac{c}{8 \sqrt{3}}$. 
This implies for the number of curves 

\begin{align*}
    m(p,q) &=(p-1)\binom{2k}{k} \\
    &\geq c' \frac{4g}{3\alpha\log(g)}\frac{g^{\frac{3}{4}\log(4)\alpha}}{\sqrt{ \alpha \log(g)}} \\
    &= \frac{4c'}{3} \frac{g^{1+\frac{3}{4}\log(4)\alpha}}{(\alpha \log(g))^{\frac{3}{2}}}, 
\end{align*}
which is strictly larger than $g^{1+\alpha}$ for large $g$, since $\frac{3}{4}\log(4) >1$.

A statistical argument provides us with a subset $\bar{\Gamma} \subset \Gamma(p,q)$ with exactly $m=[g^{1+\alpha}]$ curves, whose average intersection number does not exceed the one of $\Gamma(p,q)$. Indeed, by successively removing the curve causing the greatest number of intersection points, we obtain the system $\bar{\Gamma}$ with $\frac{\cro{\bar{\Gamma}}}{\binom{m}{2}} \leq \frac{\cro{\Gamma(p,q)}}{\binom{m(p,q)}{2}}$.
As a consequence
\begin{align*}
    \cro{\bar{\Gamma}} &\leq \frac{m(m-1)}{m(p,q)(m(p,q)-1)} \cro{\Gamma(p,q)} \\
    &\leq \frac{m(m-1)m(p,q)}{m(p,q)-1} \frac{4k}{p-1}.
\end{align*}
Since $m \leq g^{1+\alpha}$ and $m(p,q) > g^{1+a}$, for large $g$, it follows that
\begin{align*}
    \frac{m(m-1)m(p,q)}{m(p,q)-1} &\leq \frac{g^{1+\alpha}(g^{1+\alpha}-1)m(p,q)}{m(p,q)-1} \\
    &= \frac{g^{2+2\alpha}m(p,q) - g^{1+\alpha}m(p,q)}{m(p,q)-1} \\
    &< \frac{g^{2+2\alpha}m(p,q) - g^{2+2\alpha}}{m(p,q)-1}= g^{2+2\alpha}\frac{m(p,q)-1}{m(p,q)-1} \\
    &= g^{2+2\alpha}.
\end{align*} 
Finally, 
we conclude for the crossing number of $\bar{\Gamma}$:
\begin{align*}
    \cro{\bar{\Gamma}} &\leq g^{2+2\alpha} \frac{4k}{p-1} \\
    &\leq g^{2+2\alpha} \frac{3\alpha \log(g)}{\frac{4g}{3\alpha \log(g)}} \\
    &= \frac{9}{4}\alpha^2 g^{1+2\alpha} (\log(g))^2,
\end{align*}
for all but finitely many $g \in \N$.

\begin{Rem}
    All curves in $\Gamma$ are boundary curves of homeomorphic subsurfaces $\Sigma(2,k) \subset \Sigma(p,q)$. Thus, all curves are separating and topologically equivalent.
\end{Rem}

Now that we have established Theorem~\ref{polynomial}, this begs the question of whether there is true asymptotic behavior, namely:
\begin{question}
Does the limit $\displaystyle{\lim_{g \to \infty} \frac{\mathrm{Cr}(g,[g^{1+\alpha}])}{g^{1+2\alpha} (\log(g))^2}}$ exist?
\end{question}

\bigskip
\noindent
Mathematisches Institut, Universit\"at Bern, Sidlerstrasse 5, 3012 Bern, Switzerland

\smallskip
\noindent
\texttt{sebastian.baader@unibe.ch}

\smallskip
\noindent
\texttt{jasmin.joerg@unibe.ch}

\bigskip
\noindent
D\'epartement de math\'ematiques, Chemin du Mus\'ee 23, 1700 Fribourg, Switzerland

\smallskip
\noindent
\texttt{hugo.parlier@unifr.ch}


\begin{thebibliography}{99}

\bibitem{Ba}
     S.~Baader: \emph{Bipartite graphs and quasipositive surfaces}, Q.~J.~Math.~65 (2014), no.~2, 655--664.

\bibitem{BJK}
     S.~Baader, D.~Kosanovi\'c, J.~J\"org: \emph{Sparse curve systems have intermediate growth type}, arXiv:2510.14429.

\bibitem{BBS}
     S.~Baader, C.~Burrin, L.~Studer: \emph{On the crossing number of arithmetic curve systems}, to appear in Comment. Math. Helv.

\bibitem{BP}
     C.~Burrin, H.~Parlier: \emph{Modular systoles are extremal for the crossing number}, arXiv:2508.20958.

\bibitem{HP}
     A.~Hubard, H.~Parlier: \emph{Crossing number inequalities for curves on surfaces}, arXiv:2504.00916.

\bibitem{H}
     H.~Huber: \emph{Zur analytischen Theorie hyperbolischer Raumformen und Bewegungsgruppen}, Math. Ann.~138 (1959), 1--26.

\bibitem{M}
    M.~Mirzakhani: \emph{Growth of the number of simple closed geodesics on hyperbolic surfaces}, Ann. of Math. (2) 168 (2008), no.~1, 97--125.

\bibitem{PTT}
    J.~Pach, G.~Tardos, G.~Tóth: \emph{Crossings between non-homotopic edges}, J. Comb. Theory, Ser. B 156, 389--404 (2022)

\bibitem{Pa}
     H.~Parlier: \emph{Kissing numbers for surfaces}, J.~Topol.~6 (2013), no.~3, 777--791.
 
\end{thebibliography}
\end{document}